\input amstex

\documentstyle{amsppt}
\catcode`\@=11
\def\nologo{\def\logo@{}}
\catcode`\@=13 \nologo
\parindent0.7cm
\pagewidth{14 cm} \pageheight{18 cm} \pageno=1 \magnification=1200
\define\codim{\operatorname{codim}}

\topmatter

\title
On transversals of quasialgebraic families of sets
\endtitle

\author
G.R. Chelnokov, V.L. Dol'nikov
\endauthor

\affil
 Yaroslavl State University, Sovetskaya Str. 14,
Yaroslavl, 150000, Russia
\endaffil

\abstract The main results of this paper are generalizations some
classical theorems about transversals for families of finite sets to
some cases of families of infinite sets.
\endabstract
\thanks
The research of authors is supported by the Russian government
project 11.G34.31.0053 and by Russian Basic Research Foundation of
Fundamental Researches grants 10-01-00096.
\endthanks

\email grishabenruven@yandex.ru, dolnikov\@dlc.edu.yar.ru
\endemail

\endtopmatter

\document

\baselineskip=12pt

\head{ 1. Statements of the problems and the results.}
\endhead

In this paper we consider Helly-Gallai numbers for families of sets
that are similar to families of sets which are solutions for finite
systems of equations.

\definition{Definition 1}
A set $X$, $|X|\le t$, is called a {\it $t$-transversal} of a family
of sets $F$ if $A\cap X \neq \emptyset$ for every $A\in F$. By $\tau
(F)$ denote the least positive integer $t$ such that there exists a
$t$-transversal of the family $F$. This number $\tau (F)$ is called
the transversal number (or piercing number) of $P$ (see [1,10]).
\enddefinition
\definition{Definition 2}
The Helly-Gallai number $HG(t,F)$  of a family of sets $F$ is called
a minimal number $k$ such that if every subfamily $P\subseteq F$
with $|P|\le k$ has a $t$-transversal, then the family $F$ has a
$t$-transversal (see [1,2,10]).

For every family $F$, the existence of a 1-transversal is equivalent
to the condition that the intersection of all sets of~$F$ is
nonempty. Therefore a number $HG(1,F)$ is called a Helly number
$H(F)$ for a family $F$.
\enddefinition

\remark{Remark} If $F$ is a family of intervals on the line, then
$HG(t,F)=t+1$. If $F$ is a family of a convex compact sets in $\Bbb
R^d, d\ge 2$ and $t\ge 2$, then $HG(1,F)=d+1$ and numbers $HG(t,F)$
for $t \ge 2$ don't exist.
\endremark

The Helly numbers for a family of algebraic varities were found by
T.S. Motzkin [3].

\definition{Definition 3} Let $A^d_m$ be a family of sets of common zeroes in
$\Bbb R^d$ for a finite collection of polynomials of $d$ variables
and degree at most~$m$.
\enddefinition

\proclaim{Motzkin's Theorem }
$$
H(A^d_m)=\binom{m+d}{d}.
$$
\endproclaim
The Helly--Gallai numbers for algebraic varities $A^d_n$ were
determined by M.~Deza and P.~Frankl [4], and V.Dol'nikov [5]. They
are given by the formula:
$$
HG(t, A^d_m)=\binom{\binom{m+d}{d}+t-1}{t}.
$$

In the papers [5,6] the Helly--Gallai numbers for families of sets
of more general kind were considered. More precisely, families were
the zero sets of linear finite-dimensional subspaces of functions
from a ground set $V$ to in a field~$\Bbb F$.

In particular, the Helly--Gallai numbers
$$
HG(t,S_{d-1})=\binom{d+t+1}{t}
$$
for families of spheres $S_{d-1}$ in $\Bbb R^d$ were found.
Independently Helly numbers $H(S_{d-1})$ were found by H.~Maehara
[7].

Now we give some bounds for the Helly--Gallai numbers of
quasialgebraic families of sets.

\definition{Definition 4}Let $F$ be a family of sets. Denote inductively
$F^0=F$ and
$$
F^{k+1}=\{B: B=A_i\cap A_j, \,\text{ where }\, A_i,A_j \in F^k
    \,\text{ and}\, A_i \neq A_j\}.
$$
\enddefinition

\definition{Definition 5} A family of sets $F$ has the {\it
$[d,m]$-property} if $|A|\le d$ for every $A\in F^m$. If a family
$F$ has $[d,m]$-property, then such family is called a quasiagebraic
family of a dimension $m$ and a degree $d$. Denote the class of such
families by $QA^d_m$.
\enddefinition

\example {Examples}A family of lines in $\Bbb F^d$, where $\Bbb F$
is a field, or a family of lines of a finite projective plain, a
family of all edges of a simple graph $G$ or a family of sets of all
edges for a simple graph $G$ that contain a given vertex are
quasialgebraic families of a dimension $1$ and a degree $1$.

Such families $F\in QA_1^1$ are called linear families in a
literature. Let us remark that linear families investigated in
different aspects (for example see the well-known conjecture of
Erd\"os, Faber and Lovasz [8] or [9], ch. 9).

The family of circles is a quasiagebraic family of  dimension $1$
and degree $2$. The family of finite sets of cardinality $d$ is a
quasialgebraic family of  dimension $0$ and degree $d$, but also it
may be considered as a family of dimension $k$ and degree $d-k$ for
$k \ge 0$.

More generally, the families of sets, which were considered in
[5,6], are a quasialgebraic families too.

\endexample

%
%

\definition{Definition 6} A family of sets $F$ has the {\it $\{d,m\}$-property} if
$|\cap_{A\in G} A|\le d$ for every $G\subset F, |G|=m$.
\enddefinition

\remark{Remark} The family of hyperplanes in $\Bbb F^m$, where $\Bbb
F$ is an arbitrary field, has the $[1,m]$-property. The family of
hyperplanes in $\Bbb F^m$ in general position has the
$\{1,m\}$-property.
\endremark

Further just a few notes about $\{d,m\}$-property will be made and
this notion will not be used for achieving main results.

H.~Hadwiger and H.~Debrunner introduced the following concepts [1],
see also [10].

\definition{Definition 7} Let $p$ and $q$ be integers with $p\ge q\ge 2$. We
say that a family of sets $F$ has the $(p,q)$-property (in this case
we write $F\in\Pi_{p,q}$) provided $F$ has at least $p$ members and
among every $p$ members of $F$ some $q$ of them have a common point.
\enddefinition
 By $M(p,\ q;\ m)$ we denote  $\sup_F\tau(F)$ for
all finite families of convex sets $F$ in ${\Bbb R}^m$ such that
$F\in\Pi_{p,q}$. N.~Alon and D.~Kleitman [11] proved that if $m+1\le
q\le p$, then $M(p,\ q;\ m)<\infty$.

By $P(p,\ q;\ m)$ we denote  $\sup_F\tau(F)$ for all finite families
of hyperplanes in ${\Bbb R}^m$ such that $F\in\Pi_{p,q}$.

Second author proposed a conjecture [12] (see also Oberwolfach
sept.2011) that
$$
P(p,\ q;\ m)=p-q+1 \quad \text{if}\quad m+1 \le q\le p.
$$

The main purpose of this paper is to propose some general problems
of quasialgebraic families. In the paper the following theorems are
proved.

\proclaim{Theorem 1}
$$
\binom{d+m+t}{d+m}\le HG(t,QA^d_m) \le  \sum_{i=0}^{d+m}t^i.
$$
\endproclaim

\proclaim {Corollary} $H(QA^d_m)= d+m+1.$
\endproclaim

Indeed, if $t=1$, then the upper and lower estimations of Theorem 1
coincide. Thus $H(QA^d_m)= d+m+1$.

\remark{Remark} Note that the upper and lower bounds have the same
multiplicity $t^{d+m}$ by $t$.
\endremark

\proclaim{Theorem 2}
$$  \binom{d+3}{2} \le HG(2,QA^d_1)
\le  \max\left(2d^2+3,\binom{d+3}{2}\right),
$$
and
$$
\binom{d+3}{2} = HG(2,QA^d_1)
$$
for $d=1,2,3$\footnote{First author has a proof of this equality for
all $d \le 7$.}.
\endproclaim

\proclaim{Theorem 3}
$$
\binom{t+2}{2} \le HG(t,QA^1_1 ) \le
\max\left(t^2-t+3,\binom{t+2}{2}\right).
$$
\endproclaim
\remark {Remark} If $F$ is family of lines in $\Bbb F^d$, where
$\Bbb F$ is a field, or a family of lines in a projective space over
a field $\Bbb F$, then
$$
HG(t, F)=\binom{t+2}{2}\qquad [6].
$$
But if $F$ is a family of lines of an arbitrary finite projective
plain, then the authors don't know uppers estimates better than in
Theorem 3.
\endremark

\proclaim{Theorem 4} If $F \in QA^1_1$, $F \in \Pi_{p,q}$ and $|F|
\ge  \binom{p-q+3}{2}\left( \binom{p-q+2}{2}-1\right)+p-q+4$, then
$$
\tau(F) \le p-q+1.
$$
\endproclaim


\head { 2. Basic Definitions and Notations} \endhead

\definition{Definition 8} A family of sets $F$ has  {\it the  $t$-property} if a
subfamily $F\setminus \{A\}$ has a $t$-transversal $X(A)$  for every
$A\in F$ such that $X(A)\cap A=\emptyset$. Such $t$-transversal
$X(A)$ of $F\setminus \{A\}$ is called {\it a representation} of the
set $A$.
\enddefinition

\remark{Remark} Note that if $A\neq B \in F$, then $X(A)\neq X(B)$.
\endremark

\definition{Definition 9} By $a(d,m; t)$ denote the maximal $HG(t,F)$ for a family $F$
having the \linebreak $[d,m]$-property. By $b(d,m; t)$ denote the
maximal cardinality of a family $F$ having the $[d,m]$-property and
the $t$-property. By $c(d,m; t)$ denote the maximal cardinality of a
family $F$ having the $\{d,m\}$-property and the $t$-property.
\enddefinition

\remark {Remark} Later we'll prove that $a(d,m;t)$, $b(d,m;t)$ and
$c(d,m;t)$ are finite for every triple $(d,m,t)$.
\endremark

\proclaim{Proposition 1}$b(d,m; t) \ge a(d,m;t).$
\endproclaim
\demo{Proof} Assume the converse $a(d,m; t) > b(d,m; t)$. Then there
exists a family $F$ with the $(d,m)$-proper\-ty and the $t$-property
such that every $a(d,m; t)-1$ sets of $F$ have a $t$-transversal,
but $F$ doesn't. Consider a subfamily $F'\subseteq F$ of the minimal
cardinality that doesn't have a $t$-transversal. Then
$|F'|>a(d,m;)-1 \ge b(d,m; t)$. Since $F'$ is a minimal subfamily,
it follows that every $G \subset F'$ has a $t$-transversal. Thus
$F'$ has the $t$-property. Also, $F'$ has the $[d,m]$-property as a
subfamily of $F$. Hence the inequality $|F'|> b(d,m; t)$ contradicts
the definition of $b(d,m; t)$.

\enddemo

Therefore the upper bounds for $b(d,m;t)$ are also valid for $a(d,m;
t)$.

\proclaim{Proposition 2}
$$
b(d,m;t) \ge \binom{d+m+t}{t}.
$$
\endproclaim
\demo{Proof} Let $F$ be a family of all $(d+m)$-element subsets of a
$(d+m+t)$-element set. Obviously, the family $F$ has the
$[d,m]$-property and the $t$-property. Clearly, $F\setminus \{A\}$
has a $t$-transversal for all $A \in F$. Also, $F$ does not have a
$t$-transversal. It follows that
$$
a(d,m;t)\ge \binom{d+m+t}{t}.
$$
\enddemo

\definition{Definition 10} Let $F$ be a family of sets with
the $t$-property. For a set $W$ we denote
$$
G_W =\{A: A\in F, \quad W\subseteq A\}  \ \ \text{and} \ \ G_W^k
=\{A: A\in F^k, \quad W\subseteq A\}.
$$
In particular, when $W$ is a one element set $W=\{x\}$ we write
$G_x$ and $G^k_x$ respectively. Also, we denote
$$
H_x = \{ A: A\in F, \quad x\in X(A)\}.
$$
\enddefinition

\definition{Definition 11}  A set $B$ is called {\it a proper set}
for the family $F$ if there exists a non-empty subfamily $H \subset
F$ such that
$$
B = \bigcap_{A \in H} A \neq \emptyset.
$$
\enddefinition

\definition{Definition 12} Consider a set $W\neq \emptyset$. Let $B$ be a minimal proper set
such that $W \subset B$. By $\codim(W)$ denote the maximal number
$k$ such that $B \in F^{k-1}$. If such $B$ does not exists we say
that $\codim(W)=0$. If $F$ has the $[0,m]$-property, then
$\codim(W)<\infty$ of every nonempty set $W$.
\enddefinition

\head{ 3. Proof of theorem 1} \endhead

\proclaim{Lemma 1} $b(0,1;t)=t+1$.
\endproclaim
\demo{Proof} Indeed, the $[0,1]$-property of a family $F$ means that
each two sets of $F$ has an empty intersection. But $F\setminus
\{A\}$ has a $t$-transversal, then it has at most $t$ sets and $|F|
\le t+1$.

By proposition 2 $b(0,1;t)=t+1$.
\enddemo

\proclaim{Lemma 2} $b(d,1;1)= d+2$.
\endproclaim
\demo{Proof} By proposition 2 it is sufficient to prove $b(d,1;1)\le
d+2$.

Suppose a family $F$ has the $[d,1]$-property and $|F|\ge d+3$.
Consider different sets $A_1,\dots,A_{d+3}\in F$. Each
representation $X(A_i)$ consists of one element $x_i \in \cap_{j\ne
i}A_j $, and $x_i\ne x_j$ if $i\ne j$. Then
$$
A_{d+2} \cap A_{d+3} \supset \{x_1,\dots,x_{d+1}\},
$$
This contradiction completes the proof.
\enddemo

\proclaim {Lemma 3} Let $F$ be a family of sets, with the
$(d,m)$-proper\-ty and the $t$-property. Then $|G_x| \le b(d-1,m;t)$
for all $x$. Also,  $b(d,m;t) \le tb(d-1,m;t)+1$.
\endproclaim

\demo{Proof} Consider the family of sets
$$
E=\{A\setminus\{x\}: A\in G_x\}.
$$
 Let $B \in E^m$. Obviously, $B\cup\{x\} \in G_x^m
\subset F^m$. Then $|B\cup\{x\}| \le d$, therefore $|B| \le d-1$.
That means the family $E$ has the $[d,m]$-property.

Consider a set $A \in E$. Since $A\cup\{x\} \in F$, we may take
$X(A\cup\{x\})$. Then we have that $x \notin X(A\cup\{x\})$ because
$X(A\cup\{x\})\cap \left(A\cup\{x\}\right) = \emptyset$. Thus the
set $X(A\cup\{x\})$ is also a $t$-representation of the set $A$ in
the family $E$. That means the family $E$ has the $t$-property. Then
$|G_x|=|E| \le b(d,m;t)$.

Take $A \in F$ and $X(A)=\{x_1,\dots,x_t\}$. Recall that
$$
F=\{A\}\bigcup G_{x_1}\bigcup \cdots \bigcup G_{x_t},
$$
then $|F|\le1+tb(d,m;t)$.
\enddemo

Arguing as above, we see that.

\proclaim{Lemma 4} Let $F$ be a family of sets, with the
$[d,m]$-property and the $t$-property. Then
$$
|G| \le b(d-|\bigcap_{A\in G} A|,m; t) \text{\ \ for every \ } G
\subset F.
$$
\endproclaim

\proclaim{ Lemma 5} Let $F$ be a family of sets, with the
$[d,m]$-property and the $t$-property. Then
$$
|G| \le c(d-|\bigcap_{A\in G} A|,m; t) \text{\ \ for every \ } G
\subset F.
$$
\endproclaim

\proclaim{ Lemma 6} Let $B$ be a proper set. Then  $\codim(B)$
equals the maximal number $k$, for which there exist proper sets
$B_1,\dots,B_k$ such that $B=B_k\subset B_{k-1}\subset\cdots\subset
B_1$.
\endproclaim

The following inequality for $|G_x|$ is a core step in our proof of
Theorem 1.

\proclaim{Lemma 7} Let a family $F$ has the $[0,m]$-property and the
\  $t$-property. Suppose $codim(W)=k$ for a set $W$. Then
$$
|G_W| \le  t^{k-1}+\cdots+t+1.$$

In particular,
$$|G_x| \le t^{m-1}+t^{m-2}+\cdots +t+1. \ \
\text{Consequently} \ \  b(0,m;t) \le t^{m}+t^{m-1}+\cdots +t+1.$$
\endproclaim

\demo{Proof} Induction on $k$. Base for $k=1$ is obvious.

Suppose the statement is proved for all sets of the codimension $\le
k-1$. Consider a set $W$ such that $\codim(W)=k$. Assume the
converse
$$
|G_W| \ge t^{k-1}+\cdots+t+2.
$$
 Let $\{A_1,\dots,A_s\} = G_W$. Since
 $|G_x \setminus \{A_s\}| =s-1$, it follows that
there exists $y \in X(A_s)$, which belongs to at least
$$
\left\lceil \frac{s-1}{t}\right\rceil=t^{k-2}+\cdots+t+2
$$
different sets of $G_x$. Let $\bigcap_{y \in A_i}A_i=U$. Obviously,
$W \subseteq U$, $y \in U$, but $y \notin W$, thus
$$
\codim(U) \le \codim(W)-1=k-1.
$$
 This contradiction concludes the proof of the
first statement of Lemma 7.

Since $\{x\}$ isn't empty, $\codim(\{x\}) \le m$. Finally, we obtain
$$
|F| \le t|G_x|+1 \le t^{m}+t^{m-1}+\cdots +t+1.
$$
\enddemo

\proclaim{Lemma 8}
$$
b(d,m;t) \le \sum_{i=o}^{d+m} t^i.
$$
\endproclaim

\demo{Proof} By  Lemmas 3 and 7
$$
b(0,m,t) \le t|G_x|+1 \le t^{m}+t^{m-1}+\cdots +t+1.
$$
Applying recursively the inequality $b(d,m,t) \le tb(d-1,m,t)+1$
from Lemma 3  we get
$$
b(d,m,t) \le t^{m+d}+t^{m+d-1}+\cdots +t+1 =
\frac{t^{m+d+1}-1}{t-1}.
$$
\enddemo

\demo{Proof of Theorem 1} The first inequality follows from
proposition 1 and Lemma 8. The second inequality is provided by the
example in the proof of proposition 2.
\enddemo

\head  4. Proof of theorem 2
\endhead

In following three parts we will consider quasialgebraic families
for dimension 1 and degree 1 (lines); and for dimension 1, arbitrary
degree and $t=2$.

\proclaim{Lemma 9} Let $F$ be a family with the $[d,1]$-property and
the $2$-property. If
$$
|F| \ge b(d-1,1,2)+d+3,
$$
then
$$
|H_x| \le d \,\,\text{and}\,\, |F\setminus G_x|-|H_x| \le
b(d-|H_x|,1,2)
$$
for every $x$.
\endproclaim

\noindent \demo{Proof} Let $H_x =\{ A_1, \dots ,A_r\} $ and
$$
A_1=\{x,x_1\},A_2=\{x,x_2\},\dots
A_r=\{x,x_r\}.
$$

If $|H_x| \ge d+1$, then $(F \setminus G_x) \setminus \{A_1,\dots ,
A_{d+1}\}$ contains at least two sets $B$ and $C$. Therefore
$$
B \cap C \supset \{x_1, \dots , x_{d+1}\}.
$$
This contradiction
proves that $|H_x| \le d$.

Obviously
$$
\{x_1, \dots ,x_r\} \subset A \text{\ \ for every \ \ } A \in
F\setminus (G_x \cup H_x).
$$
 Applying Lemma 4 we have
$$
 |F\setminus (G_x \cup H_x)| \le b(d-|\bigcap_{A\in F\setminus (G_x \cup H_x)} A|,1,2)
\le b(d-r,1,2)=b(d-|H_x|,1,2).
$$
\enddemo

\proclaim{Lemma 10} Let $F$ be a family of sets with the
$[1,1]$-property and the $t$-property. If $|G_x|=t+1$ for some $x$,
then $x \in X(A)$ for every $A \in F\setminus G_x$ and the family
$F\setminus G_x$ has the $(t-1)$-property. Therefore $|F| \le
t+1+b(1,1,t-1)$.
\endproclaim

\demo {Proof} Assume the converse, i.e.
$$
G_x=\{A_1,\dots,A_{t+1}\},
$$
but $x \notin X(B)$ for some $B \in F\setminus G_x$. Let
$X(B)=\{x_1,\dots,x_t\}$. Since each set $A_i \bigcap X(B) \neq
\emptyset$ for $1 \le i \le t+1$, it follows that some $x_i$ belongs
to at least two sets. Without loss of generality $x_1 \in A_1\cap
A_2$. Then $\{x,x_1\}\in A_1\cap A_2$. This contradiction completes
the proof.

Thus $x \in X(A)$ for every $A \in F\setminus G_x$. It is not hard
to prove that $F\setminus G_x$ has the $(t-1)$-property. Indeed, for
every $A \in F\setminus G_x$ the family $F\setminus \left( G_x
\bigcup |{A|}\right)$ has the $(t-1)$-transversal $X(A)\setminus
\{x\}$ and $A \bigcap \left( X(A)\setminus \{x\}\right)=\emptyset$.
\enddemo

\proclaim{Lemma 11} $b(1,1,2)=6$
\endproclaim
\noindent \demo {Proof} Let $X(A)=\{x,y\}$ for some $A\in F$. By
Lemma 3 we have $|G_x| \le b(0,1,2) = 3$. If $|G_x| = 3$ for some
$x$, then by Lemma 9
$$
|F|\le 3+b(1,1,1)=6.
$$
Suppose $|G_x| < 3$ for every $x$, then

$$
|F|\le 1+ G_x +G_y \le 5.
$$
Finally, $b(1,1,2) \le 6$, hence by Proposition 2 \ $b(1,1,2)=6$.
\enddemo

\proclaim{Lemma 12} Let $F$ be a family of sets with the
$[1,1]$-property and the $t$-property. Suppose $|G_x|=|G_y|=|G_z|=t$
and $G_x,G_y,G_z$ are mutually disjoint. Then each representation
$X(A)$, \ $A \in F$, contains either all or none of elements $x,y,z$
for every $A \in F$.
\endproclaim

\noindent \demo {Proof}

Consider $x \in X(A)$ for some $A \in F$ and
$X(A)=\{x,x_1,x_2,\dots,x_{t-1}\} $. Suppose $y,z \notin X(A)$.
Since $G_y \cap G_z = \emptyset$, it follows that $A \notin G_y$ or
$A \notin G_z$. Without loss of generality $A \notin G_y$. Using
$G_x \bigcap G_y = \emptyset$, we get $x \notin B$ for every $B\in
G_y$. Then $B \cap \{x,x_1,x_2,\dots,x_{t-1}\} \neq \emptyset$ for
each $B \in G_y$. So there exist $B_1,B_2 \in G_y$ such that $B_1
\cap B_2 \cap \{x,x_1,x_2,\dots,x_{t-1}\} \neq \emptyset$. Wlog $B_1
\cap B_2 \cap \{x,x_1,x_2,\dots,x_{t-1}\} = \{x_1\}$. Then $x_1,y
\in B_1 \bigcap B_2$. This contradiction proves $y \in X(A)$.

Let $X(A) = \{x,y,x_1,\dots,x_{t-2}\}$. Since $|G_z\setminus \{A\}|
\ge t-1$ and every set of $G_z\setminus \{A\}$ contains some element
of $\{x_1,\dots,x_{t-2}\}$, it follows that two sets of
$G_z\setminus\{A\}$ contains same element. Consequently $z \in
\{x_1,\dots,x_{t-2}\}$.
\enddemo

Arguing as above we get the following statement.

\proclaim{Lemma 13} Let $F$ be a family of sets with the
$[1,1]$-property and the $t$-property. Suppose $|G_x|=|G_y|=t$ and
$G_x \cap G_y = \emptyset$. Then each representation $X(A)$, where
$A \in F \setminus (G_x \cup G_y)$, contains either all or none of
elements $x,y$.
\endproclaim

\definition{Definition 13} Let $B\neq C \in F$ and $x \in B \cap C$.
Triple $(B,C,x)$ {\it corresponds to} a representation $X(A)$ if $x
\in X(A)$.
\enddefinition

\definition {Definition 14} Consider a representation $X(A)$. By {\it the price} $P(X(A))$ we
denote the following sum
$$
P(X(A))=\sum_{x \in X(A)} \binom{|G_x|}{2}\Big/|H_x|.
$$
\enddefinition

Next two statements are obvious corollaries of definitions.

\proclaim{Lemma 14} The sum $\sum_{A \in F} P(X(A))$ equals to the
number of triples, corresponding to at least one representation.
\endproclaim

\proclaim{Lemma 15} The number of all triples is $\sum_{B \neq C \in
F}|B\cap C|$.
\endproclaim

\proclaim{Lemma 16} $b(d,1,2) \le 2d^2+3$ for $d \ge 2$.
\endproclaim

\noindent \demo {Proof}
 Assume the converse. Let $d$ be
a minimal number such that $b(d,1,2) > 2d^2+3$ and $d \ge 2$.
Consider the following cases.

First case
$$
b(d,1,2) \le b(d-1,1,2)+ d+2.
$$
Furthermore, if $d>2$, then
$$
b(d-1,1,2) \le 2(d-1)^2+3,
$$
thus
$$
b(d,1,2) \le 2(d-1)^2+3+d+2 < 2d^2+3.
$$
If $d=2$, using $b(1,1,2)=6$ from Lemma 11, we get $b(2,1,2) \le 6+
4 < 11$. The case is proven.

Assume the opposite case $b(d,1,2) \ge b(d-1,1,2)+d+3$.

Consider a family $F$ with the $[d,1]$-property and the $2$-property
of maximal cardinality. Denote $|F|=b(d,1,2)=b$. Since
$$
|B \cap C|\le d \text{\ \ for every \ \ } B,C \in F, B \neq C,
$$
it follows that
$$
 \sum_{B,C \in F,\, B \ne C}|B\cap C| \le d\frac{b(b-1)}{2}.
$$

By the other hand. Suppose $X(A)=\{x,y\}$ for some $A \in F$. From
$$
G_x \bigcup G_y = F\setminus \{A\} \quad \text{deduce} \quad
|G_x|+|G_y| \ge b-1.
$$

Also, by Lemma 9 \ $|H_x|\le d,|H_y| \le d$. Consequently
$$
P(X(A))= \binom{|G_x|}{2}\Big/|H_x| + \binom{|G_y|}{2}\Big/|H_y| \ge
\left( \binom{|G_x|}{2} + \binom{|G_y|}{2}\right)\Big/d .
$$
By the convexity of the function $\binom{x}{2}$
$$
\left( \binom{|G_x|}{2} + \binom{|G_y|}{2}\right)\Big/d \ge
2\binom{\frac{|G_x|+|G_y|}{2}}{2}/d \ge \frac{(b-1)(b-3)}{4d}.
$$

Finally, by Lemma 14 the number of all triples, corresponding to
some representation is at least $\dsize\frac{b(b-1)(b-3)}{4d}$, and
by Lemma 15 the number of all triples is at most
$\dsize\frac{b(b-1)}{2}d$. Thus we get
$$
\frac{b(b-1)(b-3)}{4d} \le \frac{b(b-1)}{2}d,
$$
that is $b \le 2d^2+3$.
\enddemo

The first part of Theorem 2 is provided by Lemma 16 and Proposition
2. For the proof of the second part we need the following lemmas.

\proclaim{Lemma 17} $b(2,1,2)=10$.
\endproclaim
\noindent \demo{Proof}  By Proposition 2 we need to prove $b(2,1,2)
\le 10$. Consider a family $F$ with the $[2,1]$-property and the
$2$-property. Suppose $|F|=11$.

By Lemmas 3 and 11 $|G_x| \le 6$ for every $x$. If $|G_x|=6$, then
$|H_x| \le 2$ by Lemma 9. If $|G_x|<6$, then $|H_x| \le 1$ by Lemma
9.

It isn't hard to see that $P(X(A)) \ge 27/2$ for every $A \in F$.
Indeed,
$$
\text{if} \ \ |G_x|=6, \ \ \text{then} \ \
\binom{|G_x|}{2}/|H_x|=15/2, \ \ \text{if} \ \ 4\le |G_y| \le 6, \ \
\text{then} \ \ \binom{|G_y|}{2}/|H_x| \ge 6.
$$

Consider a representation $X(A)=\{x,y\}$. If $|G_x|=6$, then
$$
P(X(A))=\binom{|G_x|}{2}/|H_x|+\binom{|G_y|}{2}/|H_y| \ge 15/2+6
=27/2.
$$
If both $|G_x| < 6$ and $|G_y| < 6$, then $|G_x| =|G_y| =5$ since
$|G_x|+|G_y| \ge 10$. In this case
$$
\binom{|G_x|}{2}/|H_x|+\binom{|G_y|}{2}/|H_y|=2\binom{5}{2}=20.
$$

So $ 11\times 27/2 > 2\frac{11\times 10}{2}$, we got a contradiction
with Lemmas 14 and 15.
\enddemo

Arguing as above, we get the following statement.

\proclaim{Lemma 18} $b(3,1,2)=15$
\endproclaim

\noindent \demo{Proof} By Proposition 2 we need to prove $b(3,1,2)
\le 15$. Consider a family $F$ with the $[3,1]$-property and the
$2$-property. Suppose $|F|=16$.

By Lemmas 3 and 17 $|G_x| \le 10$ for every $x$. If $|G_x|=10$, then
$|H_x| \le 3$ by Lemma 9. If $|G_x|<10$, then $|H_x| \le 2$ and if
$|G_x|<8$, then $|H_x| \le 1$ by Lemma 9. That is for values
$|G_x|=5,6,7,8,9,10$ we obtain that the values
$$
\frac{|G_x|(|G_x|-1)}{2|H_x|} \ge 10,15,21,14,18,15
$$
respectively. For any representation $\{x,y\}$ we have $|G_x|+|G_y|
\ge 15$.

We will prove
$$
\frac{|G_x|(|G_x|-1)}{2|H_x|}+\frac{|G_y|(|G_y|-1)}{2|H_y|} \ge 25.
$$
Assume the converse
$$
\frac{|G_x|(|G_x|-1)}{2|H_x|}+\frac{|G_y|(|G_y|-1)}{2|H_y|} < 25,
$$
also w.l.o.g.
$$
\frac{|G_x|(|G_x|-1)}{2|H_x|} \ge \frac{|G_y|(|G_y|-1)}{2|H_y|}.
$$
This implies
$$
\frac{|G_y|(|G_y|-1)}{2|H_y|} < \frac{25}{2},
$$
which is possible only if $|G_y|=5$. In this case $|G_x|=10$, so
$$
\frac{|G_x|(|G_x|-1)}{2|H_x|}+\frac{|G_y|(|G_y|-1)}{2|H_y|} \ge
15+10 = 25.
$$

Then for every $A \in F$
$$
P(X(A)) \ge
\frac{|G_x|(|G_x|-1)}{2|H_x|}+\frac{|G_y|(|G_y|-1)}{2|H_y|} \ge 25,
$$
which contradicts Lemmas 14 and 15.

\enddemo

Proposition 2 and Lemmas 11, 17 and 18 proves the second part of
Theorem 2.

\head{ 5. Proof of theorem 3}
\endhead

\proclaim{Lemma 19} $b(1,1,3)=10$
\endproclaim

The proof is word by word as the proof of Lemma 11.

\proclaim{Lemma 20} $b(1,1,4)=15$.
\endproclaim
\noindent\demo {Proof} By Lemma 3 $|G_x| \le 5$. If $|G_x| = 5$ for
some $x$, then by Lemma 10 we have
$$
|F| \le |G_x|+b(1,1,3)=15.
$$

If $|G_x| < 5$ for every $x$, then suppose w.l.o.g. $|F|=16$. Denote
$X(A)=\{x,y,z,t\}$ for some $A \in F$. Since
$$
F\setminus \{A\} = G_x \cup G_y \cup G_z \cup G_t \text{\ \ and all
\ \ } |G_x|,|G_y|,|G_z|,|G_t| \le 4,
$$
it follows that some three of them are mutually disjoint and have
exactly 4 elements each. Indeed, at least one of $G_x,G_y,G_z,G_t$
have 4 elements, w.l.o.g. $G_x$. Using
$$
|(G_y \setminus G_x)\cup(G_z \setminus G_x)\cup(G_t \setminus
G_x)|=11,
$$
we get that at least one of sets
$$
G_y \setminus G_x,G_z \setminus G_x,G_t
\setminus G_x
$$
have 4 elements, w.l.o.g. $G_y \setminus G_x$. Then
$$|G_y|=4 \ \ \text{and}  \ \ G_y
\cap G_x = \emptyset.
$$
Arguing analogously we get
$$
|G_z|=4 \text{\ \  and \ \ } G_z \cap G_x = G_z \cap G_y =
\emptyset.
$$

By Lemma 12 the representation of each set of $F$ contains either
all or none of elements $x,y,z$. Suppose there exists $B \ne A$ such
that $x,y,z \in X(B)$. Denote
$$
X(A)=\{x,y,z,t'\}, \quad \{C,D\}=|F\setminus (\{A,B\}\cup G_x \cup
G_y \cup G_z \cup G_t)|=2.
$$
Since $\{t,t'\}\subset C \cap D$, we got a contradiction.

Therefore we proved that $x,y,z$ belongs to the representation of
only one set of $F$. Then the price
$$
\aligned &
P(X(A))=\frac{|G_x|(|G_x-1|)}{2|H_x|}+\frac{|G_y|(|G_y-1|)}{2|H_y|}+\frac{|G_z|(|G_z-1|)}{2|H_z|}+\frac{|G_t|(|G_t-1|)}{2|H_t|}=\\&=6+6+6+\frac{|G_t|(|G_t-1|)}{2|H_t|}
>18
\endaligned
$$
for each $A \in F$. Since $18\times16 > \frac{16\times15}{2}$, we
got a contradiction with Lemmas 14 and 15.
\enddemo

 Further we use the following notation. Let the family $F$ has the
$t$-property and $A\in F$. By $x_1$ we denote $x \in X(A)$ such that
$|G_x|$ is maximal for all $x\in X(A)$. Denote by $x_i$ an element
$x \in X(A)\setminus \{x_1,\dots, x_{i-1}\}$ such that the number
$$
|G_x\setminus\big( \bigcup_{1\le j<i}G_{x_j}\big)|
$$
is maximal for all $x\in X(A)\setminus \{x_1,\dots,x_{i-1}\}$.
Denote
$$
K^A_{x_i}=G_{x_i}\setminus (G_{x_1} \cup \cdots \cup G_{x_{i-1}}) \
\  \text{and} \ \ k_i=|K^A_{x_i}|.$$

\remark {Remark } Obviously, $k_1 \ge k_2 \ge \cdots \ge k_t$. Also
if $k_i=k_j$, then $G_{x_i} \cap G_{x_j} = \emptyset$.
\endremark

By $s$ denote the maximal index such that $x_s=t$.

We need the following lemmas to prove Theorem 3.

\proclaim{Lemma 21} Let a family $F$ has the $[1,1]$-property and
the $t$-property, $|F| \ge t^2-t+4$ and $|G_x|\le t$ for every $x$.
Then
$$
\sum_{i=1}^t k_i=|F|-1,\ \ s \ge 3 \ \ \text{and} \ \ k_i \ge t -s +
2\ \ \text{for each}\ \ i.
$$
Also,
$$
\sum_{i=s}^t t-k_i \le t-3.
$$
\endproclaim

\noindent \demo{Proof} The first statement is an obvious corollary
of the definition of the numbers $k_1,\dots,k_t$.

Assume $s \le 2$. Then $k_1 \le t$, $k_2 \le t$, $k_i \le t-1$ for
$3 \le i \le t$. Thus
$$
t^2-t+3 \le |F|-1=\sum_{i=1}^t k_i \le t^2-t+2.
$$
 This contradiction
proves that $s \ge 3$.

Suppose $k_i \le t -s + 1$ for some $i$. Then $k_t \le t -s + 1$.
Recall that $k_j \le t-1$ for $s < j < t$. Thus
$$
t^2-t+3 \le |F|-1=\sum_{i=1}^t k_i \le ts +(t-1)(t-s-1)+t -s + 1\le
t^2-t+2.
$$
This contradiction proves that $k_i \ge t -s$ for each $i$.

Finally, we have
$$
\sum_{i=s}^t t-k_i  = \sum_{i=1}^t t-k_i = t^2 - |F| +1 \le t-3.
$$
\enddemo

\proclaim{Lemma 22}  Under the conditions of lemma 21, if
$X(A)\bigcap X(B)\neq \emptyset$ for some $B \in F$, then
$$
\{x_1,\dots,x_s\}\bigcap X(B)\neq \emptyset.
$$
\endproclaim
\noindent \demo{Proof} Suppose $x_i \in X(A)\bigcap X(B)$. Since
$k_i \ge t -s + 2$ by Lemma 21, we have
$$
|\left(G_{x_i}\setminus K^A_{x_i}\right)\cup \{B\}| \le s-1.
$$
Note that $K^A_{x_1}, \dots K^A_{x_s}$ are mutually disjoint, then
among the subfamilies $K^A_{x_1}, \dots K^A_{x_s}$ there exists
subfamily $K^A_{x_l}$ such that
$$
K^A_{x_l} \cap (\left(G_{x_i}\setminus K^A_{x_i}\right)\cup
\{B\})=\emptyset.
$$
 Assume that $x_l \notin X(B)$. Then each set of $K^A_{x_l}$
contains one of the elements of $X(B)$, since $B \notin K^A_{x_l}$.
Each element of $X(B)$ belongs to at most one of the sets of
$K^A_{x_l}$ and $x_i$ belongs to none of them. This contradicts the
assumption. Lemma 22 is proved.
\enddemo

\proclaim{Lemma 23}  Under the conditions of Lemma 21, if
$X(A)\bigcap X(B)\neq \emptyset$ for some $B \in F$, then
$\{x_1,\dots,x_s\}\subseteq X(B)$.
\endproclaim

\noindent \demo{Proof} Suppose $x_i \in X(A)\bigcap X(B)$. By Lemma
22 there exists some $x_l$ such that $l \le s$ and $x_l \in X(B)$.
By Lemma 21 $s \ge 3$, then Lemma 12 implies that
$\{x_1,\dots,x_s\}\subseteq X(B)$.
\enddemo

\proclaim{Lemma 24} Under the conditions of Lemma 21, \
$$
X(A) \bigcap X(B)=\emptyset \ \ \text{for every} \ \ B \in F, \ B
\neq A.$$
\endproclaim

\noindent \demo{Proof} Assume the converse. Let there exist some $i$
such that $x_i \in X(B)$. We will prove that $x_j \in X(B)$ for all
$j, \ \ 1 \le j \le t$. It means $X(A)=X(B)$, thus $A=B$.

The proof is by induction on $t-k_j$. The base claims that if there
is some $x_i \in X(B)$, then $\{x_1,\dots,x_s\} \subseteq X(B)$.
Thus the base is provided by Lemma 23.

The induction step. Consider it is proved that $x_i \in X(B)$ for
all $i$  such that $k_i \ge t-n$ for some integer nonnegative $n$.
By $q$ denote maximal $i$ such that $k_i \ge t-n$. Take some index
$j$ such that $k_j < t-n$ and $k_j$ is a maximal among all $k_l$, \,
$k_l <t -n$. Clearly, $k_j=k_{q+1}$. We have to prove $x_j \in
X(B)$.

Since by Lemma 21 \ \ $\sum_{l>q}^t t- k_l \le t-3$, it follows that
$$
t-k_j = t-k_{q+1}= \sum_{l>q}^t t- k_l -\sum_{l>q+1}^t t- k_l \le
\sum_{l>q}^t t- k_l -\sum_{l>q+1}^t 1 \le (t-3) - (t-q-1) = q - 2.
$$

Obviously, the first inequality is sharp only in two cases: $k_l
=t-1$ for $l
> q$ or $q = t -1$.

Thus $k_j \ge t -q + 2$, consequently
$$
|K^A_{x_j}\setminus\{B\}| \ge t -q +1.
$$
If $x_j \not\in X(B)$, then each set of $K^A_{x_j}\setminus\{B\}$
contains unique element of $X(B)$. Recall that $K^A_{x_i}$ are
mutually disjoint, then none of sets of $K^A_{x_j}\setminus\{B\}$
contains $x_1, \dots ,x_q$. Since there are only $t-q$ elements in
$X(B)\setminus \{x_1, \dots ,x_q\}$, we obtain the contradiction.
\enddemo

\proclaim{Lemma 25}
$$
b(1,1,t) \le \max\left(t^2-t+3,\binom{t+2}{2}\right).
$$
\endproclaim

\noindent \demo{Proof} We'll prove this by induction on $t$. The
base for $t=1,2,3,4$ is provided by Lemmas 2, 11, 19 and 20.

The step. Consider $t > 4$ and suppose the statement is proved for
all $t'$ such that $4 \le t' < t$.  Note that
$$t^2-t+3 \ge
\binom{t+2}{2} \ \  \text{for} \ \ t \ge 4.$$

Assume that there exists a family $F$ with the $[1,1]$-property and
the $t$-property such that $|F| \ge t^2-t+4$.

By Lemmas 3 and 1 $|G_{x}| \le b(0,1,t) = t+1$ for all $x$.

Assume there exists $x$  such that $G_x=t+1$. Then by Lemma 10

$$
|F| \le t+1+b(1,1,t-1) \le t+1 + (t-1)^2 - (t-1) + 3 < t^2 -t +3.
$$
Consider the opposite case. Let $|G_x| \le t$ for every $x$. Then by
Lemma 24
$$
X(A)\bigcap X(B)=\emptyset,
$$
that means $|H_x| \le 1$ for every $x$. Thus the price of any
representation
$$
\aligned & P(X(A))=\sum_{i=1}^t
\frac{|G_{x_i}|(|G_{x_i}-1|)}{2|H_{x_i}|}=\sum_{i=1}^t
\frac{|G_{x_i}|(|G_{x_i}-1|)}{2} \ge \\& \ge \sum_{i=1}^t
\frac{k_i(k_i-1)}{2} \ge t\frac12 \frac{\sum_{i=1}^t
k_i}{t}\left(\frac{\sum_{i=1}^t k_i}{t}-1 \right ) =t\frac12
\frac{|F|-1}{t}\left(\frac{|F|-1}{t}-1 \right ).
\endaligned
$$
Thus by Lemma 14 the number of triples, corresponding to some
representation is at least
$$
|F|t\frac12 \frac{|F|-1}{t}\left(\frac{|F|-1}{t}-1 \right ),
$$
 and by Lemma 15
the number of all triples is at most
$$
\frac12|F|(|F|-1).
$$

Then
$$
 |F|t\frac 12 \frac{|F|-1}{t}\left(\frac{|F|-1}{t}-1\right) \le
\frac12|F|(|F|-1),
$$
i.e. $|F| \le 2t +1$. But the inequality
$$
t^2 -t + 4 \le 2t + 1
$$
is impossible.
\enddemo

\noindent

Now Theorem 3 follows from Propositions 1 and 2 and Lemma 25.

\head{ 5. Proof of theorem 4}
\endhead

Evidently if a family $F$ has the $(p,q)$-property, then $F$ has the
$(p-1,q-1)$-property. Consequently a family $F$ has the
$(p-q+3,3)$-property. Thus it is sufficient to prove Theorem 4 when
a family $F$ has the $(p,3)$-property.

The proof will use the induction on $p$. The base for $p=3,4$.

If $p=3$ the statement is obvious. If $p=4$, then the exists $A, B,
C \in F$ such that $A\bigcap B \bigcap C = \{x\}$. We will prove
that there exists at most one $D \in F$ such that $x\notin D$.
Assume the converse. Let there exist $D\neq E \in F$ such that $x
\notin D$ and $x \notin E$. By $x_A,x_B,x_C$ denote the elements
$A\bigcap D$, $B\bigcap D$, $C\bigcap D$ respectively, if this
elements exist.

Consider sets $A,B,D,E$. Then either $A\bigcap D \bigcap E \neq
\emptyset$ or $A\bigcap D \bigcap E \neq \emptyset$ since
$$
A\bigcap B \bigcap D = A\bigcap B \bigcap E =\emptyset.
$$
Then either $x_A \in E$ or $x_B \in E$. Analogously we get that
either $x_B \in E$ or $x_C \in E$ and either $x_A \in E$ or $x_C \in
E$.

Thus $|D\bigcap E| \ge 2$. This contradiction completes the proof.

The step of the induction. Suppose $p \ge 4$. Consider a family $F$
having the $(p+1,3)$-property and
$$
|F| \ge \binom{p}{2}\left( \binom{p-1}{2}-1\right)+p+1.
$$
If the family $F$ has the $(p,3)$-property, then by the inductive
assumption $\tau(F) \le p-2 < p-1$.

Consider $Q \subset F$ such that $|Q|=p$ and $A\bigcap B \bigcap C
=\emptyset$ for every $A,B,C \in Q$. Denote
$$
W(Q)= \{x:  A\bigcap B =\{x\} \ \ \text{for} \ \ A, B \in Q\}.
$$
Then $|W(Q)|\le \binom{p}{2}$. Let $P\subset F$ and
$$
|P|= \binom{p}{2}\left( \binom{p-1}{2}-1\right)+1.
$$
Consider $A \in P$. Since the family $F$ have the
$(p+1,3)$-property, we have $A \bigcap W(Q) \neq \emptyset$. Then
there exists $x \in W(Q)$ such that
$$
|G_x|\ge 2+ \frac{|F|-p}{|W(Q)|} \ge \binom{p-1}{2}+2.
$$

Take subfamily $F\setminus G_x$. Two cases are possible:
$|F\setminus G_x| \le p-2$ and $|F\setminus G_x| \ge p-1$.

In the first case evidently, $\tau(F) \le p-1$.

In the second case we prove that the family $F\setminus G_x$ has the
$(p-1,3)$-property. Assume the converse. Then there exists $Q'
\subset F\setminus G_x$ such that $|Q'|=p-1$ and $A\bigcap B \bigcap
C =\emptyset$ for every $A,B,C \in Q'$.

Obviously $|W(Q')|\le \binom{p-1}{2}$. Then there exist $A,B \in
G_x$ such that $A\bigcap W(Q')=\emptyset$ and $B\bigcap
W(Q')=\emptyset$. Consider the family $R=\{A,B\}\bigcup Q'$. Then
$|R|=p+1$ and $A\bigcap B\bigcap C=\emptyset$ for every $A,B,C \in
R$. This contradiction proves that $F\setminus G_x$ has the
$(p-1,3)$-property.

By the inductive assumption, $\tau (F\setminus G_x) \le p-3$, then
$\tau (F) \le p-3+1=p-2 < p-1$.

\head References
\endhead

\noindent
\roster
\item"{1.}"  Hadwiger, H., Debrunner, H., Klee, V., {\it Combinatorial Geometry
in the Plane}, Holt, Rinehart and  Winston, New York, 1964.

\item"{2.}" Danzer, L., Gr\"unbaum, B., Klee, V., Helly's theorem and
its  relatives, in {\it Convexity}, pp. 101-181, Proc.~of Symposia
in Pure Math., Vol.~7,  Amer. Math. Soc., Providence, RI, 1963.

\item"{3.}"Motzkin T.S. A proof of Hilberts Nullensatz, Math. Z.
 63. (1955),  341 -- 344.

\item"{4.}"Deza M., Frankl P. A Helly type theorem for
hypersurfaces,
J. of Combin. Theory. Ser. A. (1987).

\item"{5.}"Dol'nikov V.L., A theorem of Helly type for sets defined by
    systems of equations, Mathematical Notes., V. 46, No. 5, (1989, 837 -- 840).

\item"{6.}" Dol'nikov V.L., Igonin S.A.,
    Theorems of Helly--Gallai's type (in Russian), Fundam. Prikl. Mat. 5, No.4, (1999), 1209 -- 1226 .
\item"{7.}" Maehara H.,
Helly-type theorems for spheres, Discrete\& Comput. Geom. 4, No.3,
(1989), 279 -- 285.
\item"{8.}" Erd\"os P. On the combinatorial problems I would most
to see solved, Combinatorica 1 (1981), 25 -- 42.
\item"{9.}" Aigner M., Ziegler G.M. Proofs from the Book, third edition,
Springer 2003.
\item"{10.}" Eckhoff, J., Helly, Radon, and Carath\'eodory type theorems, in
{\it Handbook of Convex Geometry}, pp. 389-448, eds. P.M. Gruber and
J.M. Wills. North-Holland, Amsterdam, 1993.
\item"{11.}" Alon, N. and Kleitman, D.J., Piercing convex sets and Hadviger --
Debrunner $(p,q)$-problem, {\it Bull. Amer. Math. Soc.}  {\bf 27}
(1992), No.  2, 252 -- 256.
\item"{12.}" Dol'nikov V.L. and Polyakova O.P., Coloring, the
$(p,q)$-property and Ramsey Theorem (in Russian), Yaroslavl State
University, Collection of surveys devoted to 30 years anniversary of
Mathematical Department, Yaroslavl, pp. 97-142, 2006.
\endroster
\enddocument

\enddocument